\documentclass{article}

\usepackage{geometry}
\usepackage{lmodern} 
\usepackage{wrapfig}
\usepackage{colortbl}
\usepackage{amsmath}
\usepackage{amssymb}
\usepackage{mathrsfs}
\usepackage{makeidx}

\usepackage{multirow}
\usepackage{caption}
\usepackage{arydshln}

\newtheorem{thm}{Theorem}[section]
\newtheorem{defin}{Definition}[section]

\usepackage{graphicx}

\usepackage{float}

\floatstyle{ruled}
\newfloat{triangle}{!ht}{tar}
\floatname{triangle}{arithmetic triangle}

\usepackage{ulem}

\newcommand\blfootnote[1]{%
  \begingroup
  \renewcommand\thefootnote{}\footnote{#1}%
  \addtocounter{footnote}{-1}%
  \endgroup
}

\usepackage{url}
\date{} 

\title{Generalized 3x + 1 Mappings : counting cycles}

\author{Robert Tremblay}

\newenvironment{changemargin}[2]{%
\begin{list}{}{%
\setlength{\leftmargin}{#1}%
\setlength{\rightmargin}{#2}%
}%
\item[]}
{\end{list}}

\begin{document}


\maketitle


\begin{abstract}
We demonstrate that the number of cycles for two problems of the family of generalized $3x + 1$ mappings is possibly finite. 
\end{abstract}

\section{Introduction}

In a previous paper~\cite{tremblay} we have determined the conditions for the existence or not of cycles for several families of generalized $3x + 1$ mappings and we have developed a method to find them. During this process there appeared a question concerning the limitation or not of the number of cycles. The answer to this question has been formulated in the form of a conjecture by many authors~\cite{matthews, carnielli} : \textsl{the number of cycles is finite}.

In this paper we study the functions that generate the infinite permutations (original Collatz problem) and $3x + 1$ problem. At first, we pick up the result we found, specifying that there can not be cycles beyond a certain value. Subsequently, we determine intrinsic properties inherent to trajectories generated by iterative application of these functions. By using these properties and the fact that these two problems are intimately linked, we will have all the necessary elements allowing us to conclude that the number of cycles produced in these two problems is limited. 



It is surely possible to carry out the search for cycles to other problems of generalized $3x + 1$ mappings as we discussed in the previous paper. We believe that the two problems that we have dealt with in details in this short paper constitute an excellent starting point in this direction.
    
\section{Original Collatz and $3x + 1$ problems}

\subsection{Functions generating these two problems~\cite{lagarias}}

Let the function $g(n)$ be defined as follows

\begin{equation}
g(n)=\left\lbrace  
\begin{array}{ll}
\frac{2n}{3} & \mbox{, if $n\equiv0\pmod{3}$}\\
\\
\frac{4n-1}{3} & \mbox{, if $n\equiv1\pmod{3}$}\\
\\
\frac{4n+1}{3} & \mbox{, if $n\equiv2\pmod{3}$}\\
\end{array}
\right.  
\label{CollatzOriginalTransfos}
\end{equation}

The iterative application of the function to integers gives rise to sequences integers, called trajectories (orbits),

\begin{equation*}
(n, g(n),g^{(2)}(n),g^{(3)}(n), \cdots, g^{(k)}(n), \cdots) ,
\end{equation*}

 with the number of iterations $k =0, 1, 2, 3, \cdots$ and $g^{(0)}(n) = n$.

The study of the iterates of g(n) is called the \textsl{the original Collatz problem}. Whe talk about \textsl{infinite permutations} because when we apply the function $g$ to all integers a first time, we find again each of them, but in a different order. Indeed, the first transformation gives the integers $2+2q$, the second $1+4q$, and the third $3+4q$, where $q$ is any integer, positive, negative or zero.



The \textsl{$3x + 1$ problem} is concerned by the iteration of the function $T(n)$, so

\begin{equation}
T(n)=\left\lbrace  
\begin{array}{ll}
\frac{n}{2} & \mbox{, if $n\equiv0\pmod{2}$}\\
\\
\frac{3n+1}{2} & \mbox{, if $n\equiv1\pmod{2}$}.\\
\end{array}
\right.  
\label{Trois_x_plus_1_Transfos}
\end{equation}
          
\subsection{Condition to the existence or not of a cycle}

A sequence of integers forms a loop when there exists a $k \ge 1$ such that 

\begin{equation}
f^{(k)}(n) = n, 
\label{condition_cycle}
\end{equation}

where $f$ represents the functions $g$ or $T$. 

If all integers in the sequence are different two by two, we have by definition a cycle of length $p = k$. Generally, we note the sequence characterizing a cycle starting with the smallest integer.

In the original Collatz problem, the first natural number forms a cycle noted $\langle1\rangle$. The following two numbers generate the cycle $\langle2, 3\rangle$ with a period $p = 2$. Two other cycles are known, namely

\begin{equation*}
\begin{array}{ccc}
\langle4, 5, 7, 9, 6\rangle & and & \langle44, 59, 79, 105, 70, 93, 62, 83, 111, 74, 99, 66\rangle,
\end{array}
\end{equation*}

 respectively, with the periods $p = 5$ and $p = 12$.

If we extend the problem from the set of natural numbers to the set of integers, we add the cycles

\begin{equation*}
\begin{array}{cccccc}
\langle0\rangle & \langle-1\rangle & \langle-2, -3\rangle & \langle-4, -5, -7, -9, -6\rangle & and & \langle-44, -59, \cdots, -66\rangle.
\end{array}
\end{equation*}

The cycles are the same with the negative integers because the function is odd, $g(-n) = -g(n)$. In addition, the cycles are closed; there are no integers other than those included in the cycles which converges towards these cycles.

In the $3x+1$ problem, for the positive integers we have the cycle $\langle1, 2\rangle$ with $p = 2$.

For the zero and negative integers we have the cycles $\langle0\rangle$, $\langle-1\rangle$, $\langle-5, -7, -10\rangle$ and the long cycle

\begin{equation*}
	\langle-17, -25, -37, -55, -82, -41, -61, -91, -136, -68, -34\rangle
\end{equation*}

with $p = 1$, $p = 3$ and $p = 11$.

The general expression giving the result of $k$ iterations of the function $f$ on an integer $n$ is 

\begin{equation}
	f^{(k)}(n) = \lambda n + \rho_k(n),  
\label{g_k_first}
\end{equation}

where

\begin{equation}
	\lambda_C = \lambda_{k_1,k_2} = \left(\frac{4}{3}\right)^{k_1}\left(\frac{2}{3}\right)^{k_2} \phantom{12} (original \phantom{1} Collatz \phantom{1} problem)  
\label{lambda_k1_k2_Collatz}
\end{equation}

and 

\begin{equation}
	\lambda_{3x+1} = \lambda_{k_3,k_4} = \left(\frac{3}{2}\right)^{k_3}\left(\frac{1}{2}\right)^{k_4} \phantom{12} (3x+1 \phantom{1} problem)   
\label{lambda_k3_k4_3x_plus_1}.
\end{equation}

In the original Collatz problem, we have 

\begin{equation}
	k_C = k_1 + k_2,
\label{k1_plus_k2_first}
\end{equation}

with $k_2$ the number of transformations of the form $2n/3$ and $k_1$, transformations of the other two kinds, $(4n\pm1)/3$.

In the $3x+1$ problem,

\begin{equation}
	k_{3x+1} = k_3 + k_4,
\label{k3_plus_k4}
\end{equation}

 with $k_3$ the number of transformations of the form $3n/2$ and $k_4$, transformations $n/2$.

Unlike parameter $\lambda$, $\rho_k(n)$ depend on the order of application of the transformations. 

Here is a brief summary of what we got in a previous paper~\cite{tremblay}, where we found the equation giving the limit condition $C$ on the smallest integer of a cycle. 


Suppose that there is a cycle of a period $p = k$, and that $m$ is its least term. Then

\begin{equation}
	m \leq \frac{par}{\frac{1}{k_i}\left|ln(\lambda)\right|} = C,
	\label{condition_RT}
\end{equation}

where $k_i = k_1$ and $par = 7/24$ in the original Collatz problem. In the $3x+1$ problem, $k_i = k_3$ and $par = 5/12$.

Essentially, the inequality~\eqref{condition_RT} specifies that the smallest integer $m$ of a cycle cannot exceed the value $C$, imposing therefore a limit on $m$. Note that $C$ increases as $\lambda$ is close to 1. Conversely, $C$ decreases very rapidly as $\lambda$ moves away from 1. 

Let $PP$ be $\lambda$ smaller than 1 ("Plus Petit que 1")  and $PG$ larger than 1 ("Plus Grand que 1"), while remaining close to 1. In writing 

\begin{equation}
	PP = 1 - \Delta PP \phantom{1234} and \phantom{1234} PG = 1 + \Delta PG,
\label{Delta_PP}
\end{equation}

we have demonstrated (theorem 2.3~\cite{tremblay}) that starting from $PP = 2/3 = 1 - 1/3$ and $PG = 4/3 = 1 + 1/3$ in the original Collatz problem, the successive products of $PP$ and $PG$ give the maxima of $C$ and gradually get closer to 1 with the increase of $k$, the total number of iterations. We have then built an algorithm that determines the conditions on $k_1$ and $k_2$ leading to the maxima of $C$. The same goes for the $3x+1$ problem. Starting from $PP = 1/2 = 1 - 1/2$ and $PG = 3/2 = 1 + 1/2$, we obtain the maxima of $C$ by carrying out the successive products of $PP$ and $PG$.

Indeed, the $\lambda$ resulting of successive products $PP\cdot PG$ is

\begin{equation}
	\lambda = PP \cdot PG = (1 + \Delta PG) \cdot (1 - \Delta PG) = 1 + \Delta PG - \Delta PP - \Delta PP\cdot \Delta PG
	\label{product}
\end{equation}	

leading to

\begin{equation}
	1-\Delta PP < 1 + \Delta PG - \Delta PP - \Delta PP\cdot \Delta PG < 1 + \Delta PG.
	\label{inequality}
\end{equation}

It is interesting to note that all $PP$ and $PG$ obtained by the successive products of $PP \cdot PG$ (except $PP = 1/2$) in the $3x + 1$ problem are the reciprocals of those obtained in the infinite permutations. $PP = 1/2$, $PG = 3/2$, $PP=3/4$, $PG=9/8$, $PP=27/32$, $PP=243/256$, $\cdots$, in the $3x + 1$ problem and $PP = 2/3$, $PG = 4/3$, $PP=8/9$, $PG=32/27$, $PG=256/243$, $\cdots$, in the problem of infinite permutations. Indeed, if we carry out the transformations 

\begin{equation*}
	k_4 \rightarrow k_1 \phantom{12} k_3 \rightarrow k_1 + k_2
\end{equation*}

in the equation~\ref{lambda_k3_k4_3x_plus_1}, we have

\begin{equation}
	\lambda_C = \frac{1}{\lambda_{3x+1}}.
	\label{reciprocals}
\end{equation}

From this property of reciprocity and from the fact that $PG$ or $PP$ is a rational number, we deduce that when 

\begin{equation*}
	\lambda_C = PG_C = \frac{N}{D}, \phantom{12} then \phantom{12} \lambda_{3x+1} = PP_{3x+1} = \frac{D}{N} \phantom{12} with \phantom{12}  N > D,  
\end{equation*} 

and

\begin{equation*}
	\Delta PG_C > \Delta PP_{3x+1}.  
\end{equation*} 

Likewise, if $\lambda_C = PP_C$, then

\begin{equation*}
	\Delta PP_C < \Delta PG_{3x+1}.  
\end{equation*}



\subsection{Periodicity}
We will show a very interesting property (hidden) resulting from the iterative application of the function $g(n)$ generating the different trajectories.

Let $n$ and $f^{(k)}(n)$ be replaced respectively by the variables $x$ and $y$ in the general expression~(\ref{g_k_first}) and using the equations~(\ref{lambda_k1_k2_Collatz}) and (\ref{lambda_k3_k4_3x_plus_1}) giving $\lambda$, 

\begin{equation}
	c = by - ax,
	\label{Diophantine_equation}
\end{equation}

In this form we have a Diophantine equation of first degree at two unknowns with

\begin{equation}
	b_C =3^{k_C}  \phantom{12} a_C = 4^{k_1} \cdot 2^{k_2}, \phantom{12} (original \phantom{1} Collatz \phantom{1} problem)  
\label{par_a_b_collatz}
\end{equation}

and 

\begin{equation}
	b_{3x+1} =2^{k_{3x+1}}  \phantom{12} a_{3x+1} = 3^{k_3}, \phantom{12} (3x+1 \phantom{1} problem)  
\label{par_a_3x_plus_1}
\end{equation}

and

\begin{equation}
	c_C = \rho_C \cdot 3^{k_C} \phantom{12} c_{3x+1} = \rho_{3x+1} \cdot 2^{k_{3x+1}}.
\label{par_c}
\end{equation}

Depending on the new parameters $a$ and $b$, the parameter $\lambda$ become

\begin{equation}
	\lambda_{a,b} = \frac{a}{b}.
\label{par_lambda}
\end{equation}

From a well-known result of Diophantine equations theory we have the theorem

\begin{thm}
Let the Diophantine equation $c = by - ax$ of first degree at two unknowns. If the coefficients $a$ and $b$ of $x$ and $y$ are prime to one another (if they have no divisor other than $1$ and $-1$ in common), this equation admits a infinity of solutions to integer values. If $(x_0, y_0)$ is a specific solution, the general solution will be $(x = x_0 + bq, y = y_0 + aq)$, where $q$ is any integer, positive, negative or zero.
	\label{equationDiophantine}
\end{thm}

\textsl{Proof}

References : Bordell\`es~\cite{bordelles}. $\blacksquare$

We may to assign to every integer of a trajectory generates by the function $T(n)$ a number $t_j = 0$ if $T^{(j)}(n)$ is even, and $t_j = 1$ if it is odd. Then, the iterative application of the function $T$ to an integer $n$ give a dyadic sequence $w_l$ of $1$ and $0$  


\begin{equation*}
w_l = (t_0, t_1, t_2, t_3, \cdots, t_j, \cdots, t_{l-1}), \phantom{1,2} with \phantom{1,2} l \ge 1.
\end{equation*} 


For a given length $l$ there are $2^l$ different dyadic sequences $w_l$ of $0$ and $1$. 

The representation of the trajectories in terms of $t_{j}$ leads to an important theorem which makes it possible to bring out an intrinsic property, namely the \textsl{periodicity}. This property has already been observed by Terras~\cite{terras} and Everett~\cite{everett} concerning the process of iterations of the function $T(n)$ generating the problem $3x+1$, and appears in a theorem which they have demonstrated by induction. We will prove it differently, using the previous theorem.


\begin{thm}
	In the $3x+1$ problem, all dyadic sequences $w_l$ of length $l = k \ge 1$ generated by any $2^l$ consecutive integers are different and are repeated periodically.
	\label{periodicity}
\end{thm}

\textsl{Proof}

Let $k = l \ge 1$ the number of iterations applied to a given integer $n$. The trajectories of length $L = l + 1$

\begin{flushleft}
	$\phantom{1,2,3,4} (T^{(0)}(n),T^{(1)}(n))$ 
	
	$\phantom{1,2,3,4} (T^{(0)}(n),T^{(1)}(n),T^{(2)}(n))$ 
	
	$\phantom{1,2,3,4} \cdots$ 
	
	$\phantom{1,2,3,4} (T^{(0)}(n),T^{(1)}(n)), \cdots, T^{(k)}(n))$
\label{trajectories}
\end{flushleft}

correspond respectively to the dyadic sequences of length $l \ge 1$ 

\begin{flushleft}
	$\phantom{1,2,3,4} w_1 = (t_0)$ 
	
	$\phantom{1,2,3,4} w_2 = (t_0, t_1)$ 
	
	$\phantom{1,2,3,4} \cdots$
	
	$\phantom{1,2,3,4} w_{l=k} = (t_0, t_1, \cdots, t_{k-1})$ .
\label{sequences}
\end{flushleft}

For a given number $l$ we have $2^l$ different dyadic sequences $w_l$ possible. 

According to theorem~\ref{equationDiophantine}, each of the $2^l$ dyadic sequences will be performed for $k = l$. Indeed, the $0$ and the $1$ of these sequences correspond to the operations on the even and odd integers. We build $2^k$ different Diophantine equations characterized by $2^k$ different combinations of the parameters $a$, $b$ and $c$, whose solutions will be given by $(x = x_0 + 2^kq, y = y_0 + m_i^{k_2}q)$. Therefore, all the integers $x_0 + 2^kq$ starting a trajectory of length $k + 1$ correspond to the same sequence $w_k$. In a sequence of $2^k$ consecutive integers, each integer must start a different sequence $w_k$, otherwise the $2^k$ different dyadic sequences will not be performed. $\blacksquare$

This theorem is interpreted as follows :

For each of integers $n$ of any $2^k$ consecutive integers we construct, from the function $T(n)$, a trajectory of length $L=k+1$ to which we associate a sequence $w_k$ of \mbox{\{0~1\}} of length $k$. The number $P$ of different sequences is exactly $P = 2^k$. Then

\begin{itemize}
	\item[-] all sequences $w_k$ appear once and only once.
	\item[-] each sequence $w_k$  is repeated periodically for any integer $n + 2^kq$ starting a trajectory, with the period $2^k$.
\end{itemize}

Given an integer $n$ and define quantities $t_{k}(n)$ by

\begin{equation}
	g^{(k)}(n)\equiv{-t_{k}(n)\pmod{3}}
\end{equation}

such that $t_{k}$ belongs to triplet of values \{-1, 0, 1\}. We could use any other triplets, for example \{0, 1, 2\}.

Then, the sequence of all integers  

\begin{equation*}
\left(
\begin{array}{m{0,5cm}m{0,3cm}m{0,3cm}m{0,3cm}m{0,3cm}m{0,3cm}m{0,3cm}m{0,3cm}m{0,3cm}m{0,3cm}m{0,3cm}m{0,3cm}m{0,3cm}m{0,5cm}}
$\cdots$ & -2 & -1 & 0 & 1 & 2 & 3 & 4 & 5 & 6 & 7 & 8  & 9 & $\cdots$ 
\end{array}
\right)
\end{equation*}

can be represented by the triadic sequence, using $n = g^{0}(n) \equiv{-t_{0}(n)\pmod{3}}$,   

\begin{equation*}
\left(
\begin{array}{m{0,5cm}m{0,3cm}m{0,3cm}m{0,3cm}m{0,3cm}m{0,3cm}m{0,3cm}m{0,3cm}m{0,3cm}m{0,3cm}m{0,3cm}m{0,3cm}m{0,3cm}m{0,5cm}}
$\cdots$ & -1 & 1 & 0 & -1 & 1 & 0 & -1 & 1 & 0 & -1 & 1 & 0 & $\cdots$ 
\end{array}
\right)
\end{equation*}

Also, each trajectory generated by iterative application of the function $g(n)$ can be represented by a triadic sequence. Then, the result of first $k$ iterations of $g(n)$ are completely described by

\begin{equation}
	w_k(n) = (t_0(n), t_1(n), \cdots, t_{k-1}(n)). 
\end{equation}

For example, the trajectories $(5, 7, 9, 6)$, $(32, 43, 57, 38)$ and all those of length $4$ starting with an integer $n = 5 + 3^3q$ where $q$ is any integer positive, negative or zero, be represented by the dyadic sequence 

\begin{equation*}
	w_3(n) = (1, -1, 0). 
\end{equation*}


\begin{thm}
	In the original Collatz problem, all sequences $w_k$ of length $k$ generated by any $3^k$ consecutive integers are different and are repeated periodically.
	\label{periodicity_Collatz}
\end{thm}

\textsl{Proof}

The proof is similar to that the theorem~\ref{periodicity} 


 $\blacksquare$     

This theorem is interpreted as follows :

For each of integers $n$ of any $3^k$ consecutive integers we construct, from the function $g(n)$, a trajectory of length $L=k+1$ to which we associate a sequence $w_k$ of \mbox{\{-1~0~1\}} of length $k$. The number $P$ of different sequences is exactly $P = 3^k$. Then

\begin{itemize}
	\item[-] all sequences $w_k$ appear once and only once.
	\item[-] each sequence $w_k$  is repeated periodically for any integer $n + 3^kq$ starting a trajectory, with the period $3^k$.
\end{itemize}

\subsection{Distribution of trajectories and average repartition }

We will determine the distribution of the trajectories of length $k$ generated by different combinations of transformations in the original Collatz and $3x+1$ problems. 


\begin{thm}
	The number of trajectories of length $k$ composed of $k_2-iterations$ of the form $2n/3$ and $k_1-iterations$ of the other two kinds, so $(4n\pm1)/3$ in the original Collatz problem, is given by

	\begin{equation}
			\eta_{k_1,k_2} = \frac{k_C!}{k_2!(k_1)!}2^{k_1}.
			\label{eta_k1_k2}
	\end{equation}	
	\label{th_eta_k1_k2}
\end{thm}

\textsl{Proof}.

The number of $k_2-combinations$ in a set with $k_C$ elements is

	\begin{equation*}
		\begin{pmatrix}
			k_C \\
			k_2
		\end{pmatrix} = \frac{k_C!}{k_2!(k_C-k_2)!} = \frac{k_C!}{k_2!k_1!}.
	\end{equation*}	

For each of these combinations we have $2^{k_1}$ combinations of $k_1$.   $\blacksquare$ 

\begin{defin}
	Defining the average repartition of the trajectories by 
	
	\begin{equation}
		R_{k_1,k_2} = \frac{3^{k_C}}{\eta_{k_1,k_2}+1}.
		\label{rep_k1_k2}
	\end{equation}
\end{defin}

In a sequence of $3^{k_C}$ consecutive integers, there are $\eta$ integers starting from the trajectories containing $k_1$ iterations of type $(4n\pm1)/3$ and $k_2$ iterations of type $2n/3$, regardless of the order of these iterations, and this $\eta$ integers are <<spaced on average>> by a value $R$. For example, let $k = 5$, $k_1 = 3$ and $k_2 = 2$. For each sequence of consecutive $243 = 3^5$ integers there are $\eta = 80$ integers whose trajectories correspond to 3 iterations of type $(4n\pm1)/3$ and 2 iterations of type $2n/3$. These integers are <<spaced on average>> by $R = 243/(80+1) = 3$. 

The number of trajectories of length $k$ composed of $k_3-iterations$ of the form $3n/2$ and $k_4-iterations$ of the form $n/2$ in the $3x+1$ problem, is given by

	\begin{equation}
			\eta_{k_3,k_4} = \frac{k_{3x+1}!}{k_3!(k_4)!}.
			\label{eta_k3_k4}
	\end{equation}

\begin{defin}
	We define the average repartition of the trajectories by 
	
	\begin{equation}
		R_{k_3,k_4} = \frac{2^{k_{3x+1}}}{\eta_{k_3,k_4}+1}.
		\label{rep_k3_k4}
	\end{equation}
\end{defin}

In the following, we will have to calculate high values of $P$, $R$ and $\eta$. For example, we can express $R$ in the original Collatz problem using natural logarithms, 


\begin{equation}
	R_{k_1,k_2} \sim k_Cln(3) - k_1ln(2) + ln(k_1!) + ln(k_2!) - ln(k_C!).
	\label{R_k1_k2_large_value}
\end{equation}

The logarithms of the factorials appearing in this last equation can be calculated by the Stirling's approximate formula, or the more accurate Ramanujan's formula

\begin{equation}
	ln(n!) \sim nln(n) - n + \frac{1}{2}ln(2\pi n),
	\label{ln_n_fact_1}
\end{equation}

\begin{equation}
	ln(n!) \sim nln(n) - n + \frac{1}{6}ln\left(8n^3 + 4n^2 + n + \frac{1}{30}\right) + \frac{1}{2}ln(\pi).
	\label{ln_n_fact_2}
\end{equation}

\subsection{Evolution of $P$, $R$ and $C$}

In the tables~\ref{NodesInfinitePermutation} and \ref{NodesInfinitePermutation_2}, we give results regarding the original Collatz problem.


In the table~\ref{NodesInfinitePermutation} we have the first values of $PP$ and $PG$ versus $k_1$ and $k_2$ giving the maxima of $C$ distributed in terms of nodes and subnodes, as presented in the previous paper. In fact, this includes the first 9 nodes. We added the natural logarithms of $C$, $R$ and $P$ as well as the exponents $r$ and $s$ in base 3 giving $\Delta PP$ and $\Delta PG$. The condition $C$ is given by the equation~\eqref{condition_RT} and the repartition $R$ (or distribution) by the equation~\eqref{rep_k1_k2}. $P = 3^k$ is the number of different trajectories for a given length $L = k + 1$, with $k = k_1 + k_2$. $P$ corresponds also to the number of consecutive integers starting the different possible sequences. 

In the table~\ref{NodesInfinitePermutation_2} we have the same information for nodes 7 to 14. We have used $\Delta PP$ and $\Delta PG$ instead of $PP$ and $PG$, by increasing the precision until the twenty-eighth decimal. This table will be useful to understand the detailed behavior of the growth of $C$ versus $P$ and $R$. 

In the tables~\ref{NodesProblem_3xPlus1} and \ref{NodesProblem_3xPlus1_2}, we give results regarding the $3x+1$ problem.

Let us now examine how the $3$ quantities $P$, $R$ and $C$ behave, one with respect to the other. Their comparative evolution should allow us to suggest an answer concerning the limitation or not of the number of cycles generated by the functions $T(n)$ and $g(n)$. In fact, as the quantities $P$ and $C$ quickly become very high, we will analyse the behavior of the logarithms of $P$, $R$ and $C$. 

\underline{Evolution of $P$}

Let $P = b^k$ (with $b = 3$ and $k = k_C$ in the original Collatz problem and $b = 2$ and $k = k_{3x+1}$ in the $3x+1$ problem), the number of different sequences associated to the trajectories generated by the functions $T(n)$ and $g(n)$. Apply the natural logarithm on each side of these equations

\begin{equation}
	ln(P_C) = k_Cln(3) \phantom{12} and \phantom{12} ln(P_{3x+1}) = k_{3x+1}ln(2).
\end{equation}

Then, \textsl{the function $ln(P)$ grows linearly with $k$}.

According to the algorithm, $\lambda_{k_1,k_2}$ (eq~\eqref{lambda_k1_k2_Collatz}) and $\lambda_{k_3,k_4}$ (eq~\eqref{lambda_k3_k4_3x_plus_1}) approach $1$ rapidly and asymptotically. 

For the original Collatz problem, we have

\begin{equation*}
	\left(\frac{4}{3}\right)^{k_1}\left(\frac{2}{3}\right)^{k_2} \sim 1,  
\end{equation*}

\begin{equation*}
	{k_1}ln(\frac{4}{3}) + {k_2}ln(\frac{2}{3}) \sim 0, \phantom{12} and  
\end{equation*}

\begin{equation*}
	\frac{k_1}{k_2} \sim -\frac{ln(\frac{2}{3})}{ln(\frac{4}{3})}.
\end{equation*}

Also,

\begin{equation*}
	k_1 + k_2 = k_C.
\end{equation*}

Resolving these last two equations,

\begin{equation}
	\frac{k_1}{k_C} \sim \frac{ln(3/2)}{ln(2)} = p_{k_1} \phantom{12} and \phantom{12} \frac{k_2}{k_C} \sim 1 - \frac{ln(3/2)}{ln(2)} = p_{k_2}.
	\label{p_k_1_and_p_k_2}
\end{equation}

These results are quickly achieved. 

For the $3x+1$ problem, we have

\begin{equation}
	\frac{k_3}{k_{3x+1}} \sim \frac{ln(2)}{ln(3)} = p_{k_3} \phantom{12} and \phantom{12} \frac{k_4}{k_{3x+1}} \sim1 - \frac{ln(2)}{ln(3)} = p_{k_4}.
\label{p_k_3_and_p_k_4}
\end{equation}

\underline{Evolution of $R$}

Now, let's analyze the growth of $R_C$ (eq~\eqref{rep_k1_k2}) in function of $k_C$ for the original Collatz Problem. 

\begin{equation*}
	R_C = \frac{3^{k_C}}{\begin{pmatrix} k_C \\ k_2 \end{pmatrix}2^{k_1} + 1} \sim \frac{3^{k_C}}{\begin{pmatrix} k_C \\ k_2 \end{pmatrix}2^{k_1}}.
\end{equation*}

Then,

\begin{equation*}
	ln(R_C) \sim k_Cln(3) - k_1ln(2) - ln\begin{pmatrix} k_C \\ k_2 \end{pmatrix},
\end{equation*}

and

\begin{equation*}
	ln(R_C) \sim k_Cln(3) - (p_{k_1}k_C)ln(2) - ln\begin{pmatrix} k_C \\ k_2 \end{pmatrix}.
\end{equation*}

The first two terms grow linearly with $k_C$. Take the last term,

\begin{equation*}
	ln\begin{pmatrix} k_C \\ k_2 \end{pmatrix} = ln\left(\frac{k_C!}{k_1!k_2!}\right) = ln(k_C!) - ln(k_1!) - ln(k_2!).
\end{equation*}

So by the Stirling's approximate formula (eq~\eqref{ln_n_fact_1})

\begin{equation*}
	ln(k!) \sim kln(k) - k + \frac{1}{2}ln(2\pi k),
\end{equation*}

we have

\begin{equation*}
	ln\begin{pmatrix} k_C \\ k_2 \end{pmatrix} \sim (-p_{k_1} ln(p_{k_1}) - p_{k_2} ln(p_{k_2})) k_C - \frac{1}{2} ln(2\pi p_{k_1} p_{k_2} k_C),
\end{equation*}

and,

\begin{equation*}
	ln(R_C) \sim k_Cln(3) - (p_{k_1}k_C)ln(2) - (-p_{k_1} ln(p_{k_1}) - p_{k_2} ln(p_{k_2})) k_C - \frac{1}{2} ln(2\pi p_{k_1} p_{k_2} k_C).
\end{equation*}

By replacing the parameters $p_{k_1}$ and $p_{k_2}$ by their respective values (equation~\ref{p_k_1_and_p_k_2}), we finally have

\begin{equation}
	ln(R_C) \sim 0.014508422k_C + \frac{1}{2} ln(1.525443029k_C).
	\label{R_C}
\end{equation} 

For the $3x+1$ problem,

\begin{equation}
	ln(R_{3x+1}) \sim 0.0346883117k_{3x+1} + \frac{1}{2} ln(1.463085787k_{3x+1}).
	\label{R_3x_plus_1}
\end{equation} 

For sufficiently high values of $k_C$ and $k_{3x+1}$, the second terms of these expressions get smaller and smaller in front of the first ones. We therefore conclude that \textsl{the logarithm of the average repartition $R$ increase linearly with the corresponding $k$}.

\underline{Evolution of $C$}

Now, let's analyze the growth of $C$ given by the equation~\eqref{condition_RT}, so 

\begin{equation*}
	m \leq k_i \cdot \frac{par}{\left|ln(\lambda)\right|} = C,
\end{equation*}

where $k_i = k_1$ and $par = 7/24$ in the original Collatz problem. In the $3x+1$ problem, $k_i = k_3$ and $par = 5/12$.

Essentially, this inequality specifies that the smallest integer $m$ of a cycle cannot exceed the value $C$, imposing therefore a limit on $m$. 

Starting from $PP = 2/3 = 1 - 1/3$ and $PG = 4/3 = 1 + 1/3$ in the original Collatz problem and from $PP = 1/2 = 1 - 1/2$ and $PG = 3/2 = 1 + 1/2$ in the $3x+1$ problem, the maxima of $C$ are given by $\lambda$ close to one and resulting of successive products $PP \cdot PG$,

\begin{equation*}
	\lambda = PP \cdot PG = (1 + \Delta PG) \cdot (1 - \Delta PG) = 1 + \Delta PG - \Delta PP - \Delta PP\cdot \Delta PG.
\end{equation*}	

Apart from the linear dependence of $C$ with the number of iterations $k_i$, the knowledge of $C$ goes through the evolution of the logarithm of the $\lambda$ versus the number of iterations.

Let us examine the expression giving $\lambda$ as a function of $\Delta PP$ and $\Delta PG$. 

Beforehand, we will prove that the successive product $PP \cdot PG$ is limited.

In a first time, we can easily see that if $\Delta PP > \Delta PG$, the product $PP \cdot PG$ will give a $\Delta PP_{new}$ such that 

\begin{equation*}
	\Delta PP_{new} > \Delta PP \cdot \Delta PG.
\end{equation*}

On the other hand, if $\Delta PP < \Delta PG$, if seems possible that the result is as small as one can imagine without being zero. In fact, if we develop $\Delta PP$ and $\Delta PG$ in base $b_C = 3$ (original Collatz problem) or in base $b_{3x+1} = 2$ ($3x+1$ problem) and pose $p \, \epsilon \, (-1,0,+1)$ we can write~\cite{tremblay}


\begin{changemargin}{2cm}{-0cm}
$\Delta PP = b^{-r} = p_{a}b^{-a} + p_{a+1}b^{-a-1} + p_{a+2}b^{-a-2} + \cdots + p_{k_{PP}}b^{-k_{PP}}$
\end{changemargin}

\begin{changemargin}{4cm}{-0cm}
$\Delta PG = b^{-s} = p_{b}b^{-b} + p_{b+1}b^{-b-1} + p_{b+2}b^{-b-2} + \cdots + p_{k_{PG}}b^{-k_{PG}}$
\end{changemargin}

with $p_a = p_b=1$ , $p_{k_{PP}} \neq 0$ and $p_{k_{PG}} \neq 0$.

In this form, the smallest possible values of $\Delta PP$ and $\Delta PG$ are respectively $b^{-k_{PP}}$ and $b^{-k_{PG}}$. Nevertheless, we will see below that the minima of $\Delta PP$ and $\Delta PG$ are higher than these possible values. The exponent, in absolute value, of the last term (the smallest) of each $\Delta PP$ (or $\Delta PG$) is equal to $k_{PP}$ (or $k_{PG}$), so the number of transformations $k=k_1+k_2$. The exponents $r$ and $s$ appear in the last columns of the tables characterizing the different nodes.


Considering the following situations, so $p$ (for precedent), $a$, $b$ and new (result of $a$ and $b$) for the Collatz and $3x+1$ problems,  

\begin{table}[h]  
\begin{center}
\begin{tabular}{cccccccc}

	\multicolumn{4}{c}{Collatz} & \multicolumn{4}{c}{3x + 1} \\
	\hline

	p &  & $\_$ &  &  &  & $\_$ &     \\ 

	  &  &  & $\cdots$      &  & $\cdots$           &  &   \\ 
	a &  &  & $\Delta PG_C$ & $>$ & $\Delta PP_{3x+1}$ &  &   \\ 

	  &  & $\cdots$      &  &  &  & $\cdots$           &   \\		
	b &  & $\Delta PP_C$ &  & $<$ &  & $\Delta PG_{3x+1}$ &   \\ 
	
	  &  &  &  &  &  &  &   \\ 
	
	new &  &  &  &  &  &  &   \\ 		
					
\end{tabular}
\end{center}
\label{EvolutionC}
\end{table}

The fact that $PP$ and $PG$ are rational and reciprocal quantities versus the two problems, allows us to write 

\begin{equation*}
	\Delta PG_C = \frac{N_a}{3^{k_{a,C}}} \phantom{2} and \phantom{2}  \Delta PP_{3x+1} = \frac{N_a}{2^{k_{a,3x+1}}},
\end{equation*} 
 
with  

\begin{equation*}
	N_a = 2^{k_{a,3x+1}} - 3^{k_{a,C}},
\end{equation*}

and $k_{a,C}$ the number of iterations in the original Collatz problem and $k_{a,3x+1}$ those in the $3x+1$ problem for the situation $a$.

Likewise,

\begin{equation*}
	\Delta PP_C = \frac{N_b}{3^{k_{b,C}}} \phantom{2} and \phantom{2} \Delta PG_{3x+1} = \frac{N_b}{2^{k_{b,3x+1}}},
\end{equation*} 
 
with  

\begin{equation*}
	N_b = 3^{k_{b,C}} - 2^{k_{b,3x+1}}.
\end{equation*}

If $\Delta PP_{3x+1} > \Delta PG_{3x+1}$, then $\Delta PG_C > \Delta PP_C$, and 

\begin{equation*}
	\Delta PG_{new,C} > \Delta PP_{new,3x+1} > \Delta PP_{3x+1} \cdot \Delta PG_{3x+1}.
\end{equation*}

If $\Delta PP_C > \Delta PG_C$, then $\Delta PG_{3x+1} > \Delta PP_{3x+1}$, and 

\begin{equation*}
	\Delta PG_{new,3x+1} > \Delta PP_{new,C} > \Delta PP_C \cdot \Delta PG_C.
\end{equation*}

Finally, if $\Delta PG_C > \Delta PP_C$, then $\Delta PG_{3x+1} > \Delta PP_{3x+1}$ is impossible.

Indeed, the first inequality leads to

\begin{equation*}
	\frac{N_a}{3^{k_{a,c}}} > \frac{N_b}{3^{k_{b,c}}}.
\end{equation*}

Like $k_b = k_a + k_p$, $k_{p,c} = (k_3)_{p,3x+1}$ and $k_3 \sim \frac{ln(2)}{ln(3)}k_{3x+1}$ for $k_{3x+1}$ sufficiently high,

\begin{equation*}
	N_b < 3^{k_{p,c}}N_a = 3^{(k_3)_{p,3x+1}}N_a \sim 3^{\frac{ln(2)}{ln(3)}k_{p,3x+1}}N_a = 2^{k_{p,3x+1}}N_a.
\end{equation*}

The second inequality leads to

\begin{equation*}
	\frac{N_b}{2^{k_{b,3x+1}}} > \frac{N_a}{2^{k_{a,3x+1}}},
\end{equation*}

and

\begin{equation*}
	N_b > 2^{k_{p,3x+1}}N_a.
\end{equation*}

which is contrary to the previous result.

We have obtained an important result which allows us to conclude that the values of $\lambda$ resulting from the successive products $PP \cdot PG$ with $PP = 1 - \Delta PP$ and $PG = 1 + \Delta PG$ are more and more close to one and the new $\Delta PP_{new}$ or $\Delta PG_{new}$ constantly decreases without ever becoming smaller than the product $\Delta PP \cdot \Delta PG$. 

We can therefore follow the evolution of $C$ versus the number of iterations $k$.

For example, if $\Delta PG >> \Delta PP$, we will have around $\Delta PG / \Delta PP$ secondary nodes and $\Delta PG$ decreases in approximate increments of $\Delta PP$. Indeed, from the equation (\ref{product}), we have

\begin{equation*}
	\Delta PG_{new} \sim \Delta PG - \Delta PP .
\end{equation*}

Let $\Delta PG_{int} = b^{-t}$ the intermediate values between $\Delta PG = b^{-s}$ and $\Delta PP = b^{-r}$; then, the exponent $t$ increases from $s$ to a value near $r$, while $k$ increases by $k_{PP}$ for each secondary node. When $\Delta PG$ approaches very close to $\Delta PP$, we have the greatest variation of the exponent $t$. Nevertheless, the new value of the exponent is never greater than the sum of $s$ and $r$. The more $\Delta PP$ close to $\Delta PG$, the more secondary nodes will follow, and the progression of the exponent $t$ in front of the number of iterations will be slowed down. It is relatively easy to be convinced of this argument by examining in detail the tables representing the primary and secondary nodes as a function of the minima of $\Delta PG$ and $\Delta PP$. 

Like $\lambda = 1 + \Delta PG$ or $\lambda = 1 - \Delta PP$, and $\Delta PG$ or $\Delta PP$ get smaller and smaller, the logarithm of $\lambda$ is approximately equal to $\Delta PG$ or $\Delta PP$ when the logarithm is developed in power series. We can then rewrite the equation~\eqref{condition_RT}, 

\begin{equation*}
	m \leq C \sim \frac{k}{b^{-t}} = k\cdot b^t ,
\end{equation*}

where $b$ is the base $2$ ($3x+1$ problem) or $3$ (original Collatz problem).

For the first $8$ nodes in the original Collatz problem, $C$ is greater than $R$. From node $9$, the values of $C$ are smaller than $R$. We will always have $C$ and $R$ smaller than $P$. Indeed, $\ln {P}$ and $\ln {R}$ grow linearly (for a sufficiently high values of $k$) and \textsl{$\ln {C}$ grows practically like a logarithm}; then, starting of node $9$, $C$ is always smaller than $R$ and the gap between the two is growing. The same thing is observed in the $3x+1$ problem, either from node $8$.

For example, for the node $N_{9,1}$ in the original Collatz problem, 

\begin{equation*}
	C = exp(12.04), \phantom{12} R = exp(17.74) \phantom{12} and \phantom{12} P = exp(1,067).
\end{equation*}

For the node $N_{14,4}$ 

\begin{equation*}
	C = exp(24.52), \phantom{12} R = exp(12,677) \phantom{12} and \phantom{12} P = exp(959,473).
\end{equation*}

For the node $N_{26,1}$ 

\begin{equation*}
	C = exp(58.25), \phantom{12} R = exp(89,401,517,209) \phantom{12} and \phantom{12} P = exp(6,770,104,587,996).
\end{equation*}


 

\subsection{Interpretation}

Take the first $b^k$ natural numbers where $k_{3x+1} = k_3 + k_4$ in the $3x+1$ problem and $k_C = k_ + k_2$ in the original Collatz problem. 

In the original Collatz problem, we have proved (theorems~\ref{periodicity_Collatz}) that the $\eta$ sequences $w$ of \{-1 0 1\} of length $k$ obtained by the transformation $g(n)$ are all present in this interval and appear only once and are repeated to all the integers $n + 3^k$, where $3^k$ is the period $P$. The $\eta$ sequences start with $\eta$ different integers. Select $k_1$ and $k_2$ in such way that $R > C$ for a $\lambda$ corresponding at a maximum value of $C$. $R$ specifies the average difference between the integers starting two consecutive sequences. We therefore expect to find very few integers between 1 and $R$ starting a sequence $w$.

The solution of the equations

\begin{equation*}
	C\left(\frac{4}{3}\right)^u \sim R \phantom{12} and \phantom{12} R\left(\frac{2}{3}\right)^v \sim C,
\end{equation*} 

makes it possible to determine the number of minimal integers $u + v$ between $m=C$ (the least integer) and $R$ being part of a cycle. All these integers start different sequences in this interval. The first equation gives the first integers of the cycle supposing that all transformations are of type $(4n \pm 1)/3$. The second equation gives the last integers supposing that all transformations are of type $2n/3$. As $R$ increases very rapidly in front of $C$, so does the number $u + v$.

For the node $N_{9,23}$ in the original Collatz problem, we have $ln(C) = 18.81$, $ln(R) = 231.37$. We get at least $u + v =1,265 $ integers in resolving the previous equations. For node $26$, this value is very large.

\textsl{We conclude that there are possibly no cycles other than the nine specified in this paper and, as the cycles are closed (that is, there are no numbers other than those belonging to cycles that end on a cycle), then the numbers such as $8, 11, 14 , \cdots$, are part of infinite trajectories}. 

The integers 8, 10, 11, 12, 13, 15, 17, 18, $\cdots$,  are in the same infinite trajectory, but the integers 14, 16, 19, $\cdots$, seem to be in other infinite trajectory.

\textsl{For the same reasons as for the original Collatz problem, we conclude that there are possibly no cycles other than the five specified in this paper. In the problem $3x + 1$ the cycles are open (for example, the number $4$ end on the cycle $\langle1, 2\rangle$). We cannot conclude that all the numbers other than those belonging to the cycles converge or not to one of five cycles}. 




\section{Conclusion}

In this paper, we have developed a method that allowed us to answer the question on the limitation or not of the number of cycles in two problems belonging to the family of generalized $3x+1$ mappings, namely the original Collatz problem (infinite permutations) and the $3x+1$ problem. We have shown that the only possible cycles are those which are already known, that is to say $9$ cycles in the first problem and $5$ in the other. 

As the function that caused the original Collatz problem generates closed cycles (there are no integers other than those included in the cycles which converges towards these cycles); then, all integers not belonging to the cycles are in infinite trajectories (divergence). In the $3x + 1$ problem, the function generates opened cycles; nevertheless, we can not be say that all integers not belonging to the cycles converge towards them, they can just as diverge. On the other hand, the natural numbers seem to converge towards the only known cycle for positive integers. In the $5x + 1$ problem, where the cycles are opened, most trajectories seem divergent. 

In conclusion, we consider that our approach opens the way to the solution of the conjecture on the limitation or not of the number of cycles for the $5x + 1$ problem and several other problems of the family of generalized $3x+1$ mappings.




\begin{flushright}
\blfootnote{e-mail:roberttremblay02@videotron.ca}
\end{flushright}



\begin{table}   
\scriptsize
\begin{center}
\begin{tabular}{|c|r|l|l|r|r|r|r|r|r|r|}
	\hline
	\multicolumn{4}{|l|}{Main nodes} & \multicolumn{7}{r|}{} \\
	\cline{1-4}
	\multirow{2}*{} & \multicolumn{3}{|l|}{Secondary nodes} & \multicolumn{7}{r|}{} \\
	\cline{2-11} 
	\ & & \multicolumn{1}{|l|}{PP} & \multicolumn{1}{|l|}{PG} & $k_1$ & $k_2$ & k & $\ln{C}$ & $\ln{R}$ & $\ln{P}$ & $r or s$ \\
	\hline
	\multicolumn{11}{|c|}{} \\
	\hline
	
	1 & 1 & \multicolumn{2}{|l|}{0.66666666666667} & 0 & 1 & 1 & & & &  \\ 
	\hline
	1 & 1 & & 1.33333333333333 & 1 & 0 & 1 & & & &  \\ 
	\hline
	\multicolumn{11}{|c|}{} \\
	\hline
	
	2 & 1 & \multicolumn{2}{|l|}{0.88888888888889} & 1 & 1 & 2 & 0.91 & 0.59 & 2.20 & 2 \\ 
	\hline		
	\multicolumn{11}{|c|}{} \\
	\hline	
	
	\multirow{2}*{3} & 1 & & 1.18518518518519 & 2 & 1 & 3 & 1.23 & 0.73 & 3.30 & 1.535 \\ 
	\cline{2-11}		
	 & 2 & & 1.05349794238683 & 3 & 2 & 5 & 2.82 & 1.10 & 5.49 & 2.665 \\ 
	\hline			
	\multicolumn{11}{|c|}{} \\
	\hline
	
	\multirow{2}*{4} & 1 & \multicolumn{2}{|l|}{0.93644261545496} & 4 & 3 & 7 & 2.88 & 1.36 & 7.69 & 2.508 \\ 
	\cline{2-11}		
	 & 2 & \multicolumn{2}{|l|}{0.98654036854514} & 7 & 5 & 12 & 5.02 & 1.66 & 13.18 & 3.921 \\ 
	\hline			
	\multicolumn{11}{|c|}{} \\
	\hline			
	
	\multirow{3}*{5} & 1 & & 1.03931824834386 & 10 & 7 & 17 & 4.33 & 1.87 & 18.68 & 2.946 \\ 
	\cline{2-11}		
	 & 2 & & 1.02532940775684 & 17 & 12 & 29 & 5.29 & 2.31 & 31.86 & 3.346 \\ 
	\cline{2-11}			
	 & 3 & & 1.01152885180861 & 24 & 17 & 41 & 6.41 & 2.66 & 45.04 & 4.062 \\ 
	\hline		
	\multicolumn{11}{|c|}{} \\
	\hline

	6 & 1 & \multicolumn{2}{|l|}{0.99791404625731} & 31 & 22 & 53 & 8.37 & 2.97 & 58.23 & 5.618 \\ 
	\hline		
	\multicolumn{11}{|c|}{} \\
	\hline

	\multirow{5}*{7} & 1 & & 1.00941884941434 & 55 & 39 & 94 & 7.44 & 3.84 & 103.27 & 4.246 \\ 
	\cline{2-11}		
	 & 2 & & 1.00731324838746 & 86 & 61 & 147 & 8.14 & 4.84 & 161.50 & 4.477 \\ 
	\cline{2-11}			
	 & 3 & & 1.00521203954693 & 117 & 83 & 200 & 8.79 & 5.76 & 219.72 & 4.785 \\ 
	\cline{2-11}			
	 & 4 & & 1.00311521373084 & 148 & 105 & 253 & 9.54 & 6.65 & 277.95 & 5.253 \\ 
	\cline{2-11}			
	 & 5 & & 1.00102276179641 & 179 & 127 & 306 & 10.84 & 7.51 & 336.18 & 6.267 \\ 
	\hline		
	\multicolumn{11}{|c|}{} \\
	\hline
	
	\multirow{2}*{8} & 1 & \multicolumn{2}{|l|}{0.99893467461992} & 210 & 149 & 359 & 10.96 & 8.36 & 394.40 & 6.230 \\ 
	\cline{2-11}		
	 & 2 & \multicolumn{2}{|l|}{0.99995634684222} & 389 & 279 & 665 & 14.77 & 13.11 & 730.58 & 9.138 \\ 
	\hline			
	\multicolumn{11}{|c|}{} \\
	\hline

	\multirow{5}*{9} & 1 & & 1.00097906399185 & 568 & 403 & 971 & 12.04 & 17.74 & 1,066.75 & 6.307 \\ 
	\cline{2-11}		
	 & 2 & & 1.00093536809484 & 957 & 679 & 1,636 & 12.61 & 27.65 & 1,797.33 & 6.349 \\ 
	\cline{2-11}			
	 & \dots & &  &  &  &  &  &  &   \\ 
	\cline{2-11}			
	 & 22 & & 1.00006185061131 & 8,737 & 6,199 & 14,936 & 17.53 & 221.70 & 16,408.87 & 8.821 \\ 
	\cline{2-11}			
	 & 23 & & 1.00001819475356 & 9,126 & 6,475 & 15,601 & 18.81 & 231.37 & 17,139.45 & 9.935 \\ 
	\hline		
				
\end{tabular}
\end{center}
\caption{Nodes - Infinite Permutations - Nodes 1 to 9}
\label{NodesInfinitePermutation}
\normalsize
\end{table}

\begin{table}   
\scriptsize
\begin{center}
\rotatebox{90}{
\begin{tabular}{|c|r|l|l|r|r|r|r|r|r|r|}
	\hline
	\multicolumn{4}{|l|}{Main nodes} & \multicolumn{7}{r|}{} \\
	\cline{1-4}
	\multirow{2}*{} & \multicolumn{3}{|l|}{Secondary nodes} & \multicolumn{7}{r|}{} \\
	\cline{2-11} 
	\ & & \multicolumn{1}{|l|}{$\Delta PP$} & \multicolumn{1}{|l|}{$\Delta PG$} & $k_1$ & $k_2$ & k & $\ln{C}$ & $\ln{R}$ & $\ln{P}$ & $r or s$ \\
	\hline
	\multicolumn{11}{|c|}{} \\
	\hline
	\multirow{3}*{7} & $\cdots$ & &  &  &  & & & &  \\ 
	\cline{2-11}
	  & 4 &  & 0.0031152137308416658467349706 & 148 & 105 & 253 & 9.5381 & 6.647 & 277.9 & 5.253407026 \\ 
	\cline{2-11}
	  & 5 &  & 0.0010227617964117672208313996 & 179 & 127 & 306 & 10.8410 & 7.512 & 336.2 & 6.267223422 \\  
	\hline
	\multicolumn{11}{|c|}{} \\
	\hline	
	\multirow{2}*{8} & 1 & \multicolumn{2}{|l|}{0.0010653253800741109929206204} & 210 & 149 & 359 & 10.9589 & 8.362 & 394.4 & 6.230109635 \\ 
	\cline{2-11}		
	 & 2 & \multicolumn{2}{|l|}{0.0000436531577618341853224779} & 389 & 276 & 665 & 14.7706 & 13.109 & 730.6 & 9.138105444 \\ 
	\hline	
	\multicolumn{11}{|c|}{} \\
	\hline

	\multirow{23}*{9} & 1 & & 0.0009790639918678842653176021 & 568 & 403 & 971 & 12.0394 & 17.737 & 1,066.8 & 6.306968914 \\ 
	\cline{2-11}		
	 & 2 & & 0.0009353680948711569096002737 & 957 & 679 & 1,636 & 12.6067 & 27.645 & 1,797.3 & 6.348527587 \\ 
	\cline{2-11}	
	 & 3 & & 0.0008916741053383146967578837 & 1,346 & 955 & 2,301 & 12.9956 & 37.463 & 1,797.3 & 6.392072906 \\ 
	\cline{2-11}		
	 & 4 & & 0.0008479820231860908048872943 & 1,735 & 1,231 & 2,966 & 13.2997 & 47.238 & 3,258.5 & 6.437804487 \\ 
	\cline{2-11}		
	 & 5 & & 0.0008042918483312220469450805 & 2,124 & 1,507 & 3,631 & 13.5549 & 56.986 & 3,989.1 & 6.485953628 \\ 
	\cline{2-11}		
	 & 6 & & 0.0007606035806904488705888572 & 2,513 & 1,783 & 4,296 & 13.7789 & 66.718 & 4,719.6 & 6.536790394 \\ 
	\cline{2-11}	
	 & 7 & & 0.0007169172201805153580186127 & 2,902 & 2,059 & 4,961 & 13.9819 & 76.437 & 5,450.2 & 6.590632794 \\ 
	\cline{2-11}		
	 & 8 & & 0.0006732327667181692258180498 & 3,291 & 2,335 & 5,626 & 14.1706 & 86.148 & 6,180.8 & 6.647858848 \\ 
	\cline{2-11}	
	 & 9 & & 0.0006295502202201618247959332 & 3,680 & 2,611 & 6,291 & 14.3493 & 95.851 & 6,911.4 & 6.708922707 \\ 
	\cline{2-11}		
	 & 10 & & 0.0005858695806032481398274443 & 4,069 & 2,887 & 6,956 & 14.5217 & 105.549 & 7,641.9 & 6.774376574 \\ 
	\cline{2-11}	
	 & 11 & & 0.0005421908477841867896955426 & 4,458 & 3,163 & 7,621 & 14.6905 & 115.242 & 8,372.5 & 6.844901134 \\ 
	\cline{2-11}		
	 & 12 & & 0.0004985140216797400269323340 & 4,847 & 3,439 & 8,286 & 14.8581 & 124.932 & 9,103.1 & 6.921348796 \\ 
	\cline{2-11}		
	 & 13 & & 0.0004548391022066737376604466 & 5,236 & 3,715 & 8,951 & 15.0270 & 134.618 & 9,833.7 & 7.004806793 \\ 
	\cline{2-11}		
	 & 14 & & 0.0004111660892817574414344127 & 5,625 & 3,991 & 9,616 & 15.1996 & 144.301 & 10,564.3 & 7.096692250 \\ 
	\cline{2-11}	
	 & 15 & & 0.0003674949828217642910820580 & 6,014 & 4,267 & 10,281 & 15.3787 & 153.982 & 11,294.8 & 7.198900807 \\ 
	\cline{2-11}		
	 & 16 & & 0.0003238257827434710725458978 & 6,403 & 4,543 & 10,946 & 15.5678 & 163.661 & 12,025.4 & 7.314049713 \\ 
	\cline{2-11}				
	 & 17 & & 0.0002801584889636582047245402 & 6,792 & 4,819 & 11,611 & 15.7717 & 173.338 & 12,756.0 & 7.445898036 \\ 
	\cline{2-11}	
	 & 18 & & 0.0002364931013991097393140960 & 7,181 & 5,095 & 12,276 & 15.9968 & 183.013 & 13,486.6 & 7.600125728 \\ 
	\cline{2-11}		
	 & 19 & & 0.0001928296199666133606495956 & 7,570 & 5,371 & 12,941 & 16.2536 & 192.687 & 14,217.1 & 7.785916511 \\ 
	\cline{2-11}		
	 & 20 & & 0.0001491680445829603855464127 & 7,959 & 5,647 & 13,606 & 16.5604 & 202.360 & 14,947.7 & 8.019605426 \\ 
	\cline{2-11}		
	 & 21 & & 0.0001055083751649457631416954 & 8,348 & 5,923 & 14,271 & 16.9544 & 212.031 & 15,678.3 & 8.334805934 \\ 
	\cline{2-11}	
	 & 22 & & 0.0000618506116293680747358036 & 8,737 & 6,199 & 14,936 & 17.5340 & 221.702 & 16,408.9 & 8.820935894 \\ 
	\cline{2-11}		
	 & 23 & & 0.0000181947538930295336337538 & 9,126 & 6,475 & 15,601 & 18.8011 & 231.371 & 17,139.5 & 9.934694310 \\ 
	\hline
	\multicolumn{11}{|c|}{} \\	
	\hline

	\multirow{2}*{10} & 1 & \multicolumn{2}{|l|}{0.00002545911981272640150133296} & 9,515 & 6,751 & 16,266 & 18.5069 & 241.039 & 17870.0 & 9.628905092 \\ 
	\cline{2-11}
	 & 2 & \multicolumn{2}{|l|}{0.00000726490745807872082267250} & 18,641 & 13,226 & 31,867 & 20.4334 & 467.708 & 35,009.5 & 10.77036469  \\  
	\hline

	\multicolumn{11}{|c|}{} \\	
	\hline

	\multirow{2}*{11} & 1 & & 0.0000109297142517475574299296 & 27,767 & 19,701 & 47,468 & 20.4235 & 694.239 & 52,148.9 & 10.39859604 \\ 
	\cline{2-11}		
	 & 2 & & 0.000003664727390306254413089 & 46,408 & 32,927 & 79,335 & 22.0298 & 1,156.808 & 87,158.4 & 11.39324285 \\ 
	\hline

	\multicolumn{11}{|c|}{} \\	
	\hline

	12 & 1 & \multicolumn{2}{|l|}{0.0000036002066916778116074911} & 65,049 & 46,153 & 111,202 & 22.3853 & 1,619.289 & 122,167.9 & 11.40941115 \\ 
	\hline

	\multicolumn{11}{|c|}{} \\	
	\hline

	13 & 1 & & 0.0000000645075048523645826212 & 111,457 & 79,080 & 190,537 & 26.9457 & 2,770.514 & 209,326.3 & 15.07036143 \\ 
	\hline

	\multicolumn{11}{|c|}{} \\	
	\hline

	\multirow{6}*{14} & 1 & \multicolumn{2}{|l|}{0.000003535699419065802125212} & 176,506 & 125,233 & 301,739 & 23.4016 & 4,384.012 & 331,494.2 & 11.42586839 \\ 
	\cline{2-11}
	 & 2 & \multicolumn{2}{|l|}{0.0000034711921422925849744807} & 287,963 & 204,313 & 492,276 & 23.9095 & 7,148.481 & 540,820.5 & 11.44262867  \\  
	\cline{2-11}
	 & 3 & \multicolumn{2}{|l|}{0.0000034066848613581643542882} & 399,420 & 283,393 & 682,813 & 24.2554 & 9,912.869 & 750,146.8 & 11.45970335  \\ 
	\cline{2-11}
	 & 4 & \multicolumn{2}{|l|}{0.0000033421757626253399962058} & 510,877 & 362,473 & 873,350 & 24.5206 & 12,677.217 & 959,473.0 & 11.47710446  \\ 
	\cline{2-11}
	& $\cdots$ &  \multicolumn{2}{|l|}{}  &  &  &  &  &  &  & \\  
	\cline{2-11}
	 & 55 & \multicolumn{2}{|l|}{0.0000000523005186232530720965} & $\cdots$ & $\cdots$ & $\cdots$ & 31.1734 & 153,654 & 11,635,114 & 15.26130713\\ 
	\hline
					
\end{tabular}
}
\end{center}
\caption{Nodes - Infinite Permutations - Nodes 7 to 14}
\label{NodesInfinitePermutation_2}
\normalsize
\end{table}


\begin{table}  
\begin{center}
\begin{tabular}{|c|r|l|l|r|r|r|r|r|r|}
	\hline
	\multicolumn{4}{|l|}{Main nodes} & \multicolumn{6}{r|}{} \\
	\cline{1-4}
	\multirow{2}*{} & \multicolumn{3}{|l|}{Secondary nodes} & \multicolumn{6}{r|}{} \\
	\cline{2-10} 
	\ & & \multicolumn{1}{|l|}{PP} & \multicolumn{1}{|l|}{PG} & $k_3$ & $k_4$ & k & $\ln{C}$ & $\ln{R}$ & $\ln{P}$ \\
	\hline
	\multicolumn{10}{|c|}{} \\
	\hline
	
	1 & 1 & \multicolumn{2}{|l|}{0.50000000000000} & 0 & 1 & 1 & & &   \\ 
	\hline
	1 & 1 & & 1.500000000000000 & 1 & 0 & 1 & & &   \\ 
	\hline
	\multicolumn{10}{|c|}{} \\
	\hline
	
	2 & 1 & \multicolumn{2}{|l|}{0.75000000000000} & 1 & 1 & 2 & 0.37 & 0.29 & 1.39  \\ 
	\hline		
	\multicolumn{10}{|c|}{} \\
	\hline	
	
	3 & 1 & & 1.12500000000000 & 2 & 1 & 3 & 1.96 & 0.69 & 2.08  \\ 
	\hline		
	\multicolumn{10}{|c|}{} \\
	\hline
	
	\multirow{2}*{4} & 1 & \multicolumn{2}{|l|}{0.84375000000000} & 3 & 2 & 5 & 2.00 & 1.07 & 3.47  \\ 
	\cline{2-10}		
	 & 2 & \multicolumn{2}{|l|}{0.94921875000000} & 5 & 3 & 8 & 3.69 & 1.50 & 5.55  \\ 
	\hline			
	\multicolumn{10}{|c|}{} \\
	\hline			

	\multirow{2}*{5} & 1 & & 1.06787109375000 & 7 & 4 & 11 & 3.79 & 1.82 & 7.62  \\ 
	\cline{2-10}		
	 & 2 & & 1.01364326477050 & 12 & 7 & 19 & 5.91 & 2.34 & 13.17  \\ 
	\hline			
	\multicolumn{10}{|c|}{} \\
	\hline	
	
	\multirow{3}*{6} & 1 & \multicolumn{2}{|l|}{0.96216919273138} & 17 & 10 & 27 & 5.21 & 2.77 & 18.72  \\ 
	\cline{2-10}		
	 & 2 &\multicolumn{2}{|l|}{0.97529632178184} & 29 & 17 & 46 & 6.18 & 3.69 & 31.88  \\ 
	\cline{2-10}			
	 & 3 & \multicolumn{2}{|l|}{0.98860254772961} & 41 & 24 & 65 & 7.31 & 4.53 & 45.05  \\ 
	\hline		
	\multicolumn{10}{|c|}{} \\
	\hline

	7 & 1 & & 1.00209031404109 & 53 & 31 & 84 & 9.27 & 5.32 & 58.22  \\ 
	\hline		
	\multicolumn{10}{|c|}{} \\
	\hline

	\multirow{5}*{8} & 1 & \multicolumn{2}{|l|}{0.99066903751619} & 94 & 55 & 149 & 8.34 & 7.86 & 103.28  \\ 
	\cline{2-10}		
	 & 2 & \multicolumn{2}{|l|}{0.99273984691538} & 147 & 86 & 233 & 9.04 & 10.99 & 161.50  \\ 
	\cline{2-10}			
	 & 3 & \multicolumn{2}{|l|}{0.99481498495653} & 200 & 117 & 317 & 9.68 & 14.06 & 219.73  \\ 
	\cline{2-10}			
	 & 4 & \multicolumn{2}{|l|}{0.99689446068787} & 253 & 148 & 401 & 10.43 & 17.10 & 277.95  \\ 
	\cline{2-10}			
	 & 5 & \multicolumn{2}{|l|}{0.99897828317652} & 306 & 179 & 485 & 11.73 & 20.11 & 336.18  \\ 
	\hline		
	\multicolumn{10}{|c|}{} \\
	\hline
	
	\multirow{2}*{9} & 1 & & 1.00106646150859 & 359 & 210 & 569 & 11.85 & 23.10 & 394.40  \\ 
	\cline{2-10}		
	 & 2 & & 1.00004365506344 & 665 & 389 & 1,054 & 15.66 & 40.23 & 730.58  \\ 
	\hline			
	\multicolumn{10}{|c|}{} \\
	\hline

	\multirow{5}*{10} & 1 & \multicolumn{2}{|l|}{0.99902189363685} & 971 & 568 & 1,539 & 12.93 & 57.24 & 1,066.75  \\ 
	\cline{2-10}		
	 & 2 & \multicolumn{2}{|l|}{0.99906550600100} & 1,636 & 957 & 2,593 & 13.50 & 94.07 & 1,797.33  \\ 
	\cline{2-10}			
	 & \dots & &  &  &  &  &  &  &   \\ 
	\cline{2-10}			
	 & 22 & \multicolumn{2}{|l|}{0.99993815321363} & 14,936 & 8,737 & 23,673 & 18.43 & 826.40 & 16,408.87  \\ 
	\cline{2-10}			
	 & 23 & \multicolumn{2}{|l|}{0.99998180557715} & 15,601 & 9,126 & 24,727 & 19.69 & 862.98 & 17,139.45  \\ 
	\hline
				
\end{tabular}
\end{center}
\caption{Nodes - 3x + 1}
\label{NodesProblem_3xPlus1}
\end{table}

\begin{table}   
\scriptsize
\begin{center}
\rotatebox{90}{
\begin{tabular}{|c|r|l|l|r|r|r|r|r|r|r|}
	\hline
	\multicolumn{4}{|l|}{Main nodes} & \multicolumn{7}{r|}{} \\
	\cline{1-4}
	\multirow{2}*{} & \multicolumn{3}{|l|}{Secondary nodes} & \multicolumn{7}{r|}{} \\
	\cline{2-11} 
	\ & & \multicolumn{1}{|l|}{$\Delta PP$} & \multicolumn{1}{|l|}{$\Delta PG$} & $k_3$ & $k_4$ & k & $\ln{C}$ & $\ln{R}$ & $\ln{P}$ & $r or s$ \\
	\hline
	\multicolumn{11}{|c|}{} \\
	\hline
	\multirow{3}*{8} & $\cdots$ & &  &  &  & & & &  \\ 
	\cline{2-11}
	  & 4 & \multicolumn{2}{|l|}{0.0031055393121348348272949815} & 253 & 148 & 401 & 10.4309 & 17.10 & 277.95 & 8.330940454 \\ 
	\cline{2-11}
	  & 5 & \multicolumn{2}{|l|}{0.0010217168234779627751601743} & 306 & 179 & 485 & 11.7339 & 20.11 & 336.18 & 9.934788887 \\  
	\hline
	\multicolumn{11}{|c|}{} \\
	\hline	
	\multirow{2}*{9} & 1 &  & 0.0010664615085860798682402781  & 359 & 210 & 569 & 11.8518 & 23.10 & 394.40 & 9.872952389 \\ 
	\cline{2-11}		
	  & 2 & & 0.0000436550634432030074328558 & 665 & 389 & 1,054 & 15.6635 & 40.23 & 730.58 & 14.48349148 \\ 
	\hline	
	\multicolumn{11}{|c|}{} \\
	\hline

	\multirow{23}*{10} & 1 & \multicolumn{2}{|l|}{0.0009781063631475096860402899} & 971 & 568 & 1,539 & 12.9323 & 57.24 & 1,066.75 & 9.997721021 \\ 
	\cline{2-11}		
	 & 2 & \multicolumn{2}{|l|}{0.0009344939989996440837028075} & 1,636 & 957 & 2,593 & 13.4996 & 94.07 & 1,797.33 & 10.06352698 \\ 
	\cline{2-11}	
	 & 3 & \multicolumn{2}{|l|}{0.0008908797309512546982201879} & 2,301 & 1,346 & 3,647 & 13.8885 & 130.80 & 2,527.91 & 10.13248170 \\ 
	\cline{2-11}		
	 & 4 & \multicolumn{2}{|l|}{0.0008472635589192266314371305} & 2,966 & 1,735 & 4,701 & 14.1926 & 167.49 & 3,258.48 & 10.20490156 \\ 
	\cline{2-11}		
	 & 5 & \multicolumn{2}{|l|}{0.0008036454828204413568121826} & 3,631 & 2,124 & 5,755 & 14.4477 & 204.15 & 3,989.06 & 10.28115316 \\ 
	\cline{2-11}		
	 & 6 & \multicolumn{2}{|l|}{0.0007600255025717767192593421} & 4,296 & 2,513 & 6,809 & 14.6717 & 240.79 & 4,719.64 & 10.36166455 \\ 
	\cline{2-11}	
	 & 7 & \multicolumn{2}{|l|}{0.0007164036180901069349896530} & 4,961 & 2,902 & 7,863 & 14.8748 & 277.43 & 5,450.22 & 10.44693976 \\ 
	\cline{2-11}		
	 & 8 & \multicolumn{2}{|l|}{0.0006727798292923025913527942} & 5,626 & 3,291 & 8,917 & 15.0634 & 314.05 & 6,180.79 & 10.53757793 \\ 
	\cline{2-11}	
	 & 9 & \multicolumn{2}{|l|}{0.0006291541360952306466786615} & 6,291 & 3,680 & 9,971 & 15.2422 & 350.67 & 6,911.37 & 10.63429887 \\ 
	\cline{2-11}		
	 & 10 & \multicolumn{2}{|l|}{0.0005855265384157544301189422} & 6,956 & 4,069 & 11,025 & 15.4146 & 387.28 & 7,641.95 & 10.73797782 \\ 
	\cline{2-11}	
	 & 11 & \multicolumn{2}{|l|}{0.0005418970361707336414886834} & 7,621 & 4,458 & 12,079 & 15.5834 & 423.89 & 8,372.52 & 10.84969362 \\ 
	\cline{2-11}		
	 & 12 & \multicolumn{2}{|l|}{0.0004982656292770243511078528} & 8,286 & 4,847 & 13,133 & 15.7510 & 460.49 & 9,103.10 & 10.97079732 \\ 
	\cline{2-11}		
	 & 13 & \multicolumn{2}{|l|}{0.0004546323176514789996428931} & 8,951 & 5,236 & 14,187 & 15.9198 & 497.09 & 9,833.68 & 11.10301214 \\ 
	\cline{2-11}		
	 & 14 & \multicolumn{2}{|l|}{0.0004109971012109463979482692} & 9,616 & 5,625 & 15,241 & 16.0924 & 533.69 & 10,564.26 & 11.24858416 \\ 
	\cline{2-11}	
	 & 15 & \multicolumn{2}{|l|}{0.0003673599798722717269080083} & 10,281 & 6,014 & 16,295 & 16.2716 & 570.28 & 11,294.83 & 11.41051791 \\ 
	\cline{2-11}		
	 & 16 & \multicolumn{2}{|l|}{0.0003237209535522965372772337} & 10,946 & 6,403 & 17,349 & 16.4607 & 606.88 & 12,025.41 & 11.59296163 \\ 
	\cline{2-11}				
	 & 17 & \multicolumn{2}{|l|}{0.0002800800221678587495236908} & 11,611 & 6,792 & 18,403 & 16.6645 & 643.47 & 12,755.99 & 11.80187330 \\ 
	\cline{2-11}	
	 & 18 & \multicolumn{2}{|l|}{0.0002364371856357926536692672} & 12,276 & 7,181 & 19,457 & 16.8896 & 680.06 & 13,486.56 & 12.04625543 \\ 
	\cline{2-11}		
	 & 19 & \multicolumn{2}{|l|}{0.0001927924438729289091315050} & 12,941 & 7,570 & 20,511 & 17.1465 & 716.64 & 14,217.14 & 12.34066387 \\ 
	\cline{2-11}		
	 & 20 & \multicolumn{2}{|l|}{0.0001491457967960945445651067} & 13,606 & 7,959 & 21,565 & 17.4533 & 753.23 & 14,947.72 & 12.71098906 \\ 
	\cline{2-11}		
	 & 21 & \multicolumn{2}{|l|}{0.0001054972443221129577034338} & 14,271 & 8,348 & 22,619 & 17.8473 & 789.82 & 15,678.30 & 13.21050706 \\ 
	\cline{2-11}	
	 & 22 & \multicolumn{2}{|l|}{0.0000618467863678039151999991} & 14,936 & 8,737 & 23,673 & 18.4269 & 826.40 & 16,408.87 & 13.98094184 \\ 
	\cline{2-11}		
	 & 23 & \multicolumn{2}{|l|}{0.0000181944228499835524699513} & 15,601 & 9,126 & 24,727 & 19.6940 & 862.98 & 17,139.45 & 15.74614419 \\ 
	\hline
	\multicolumn{11}{|c|}{} \\	
	\hline

	\multirow{2}*{11} & 1 & & 0.0000254598463145356264684468 & 16,266 & 9,515 & 25,781 & 19.3998 & 899.57 & 17,870.03 & 15.26141676 \\ 
	\cline{2-11}
	 & 2 & & 0.0000072649602373425317419577 & 31,867 & 18,641 & 50,508 & 21.3263 & 1,757.64 & 35,009.48 & 17.07061367  \\  
	\hline

	\multicolumn{11}{|c|}{} \\	
	\hline

	\multirow{2}*{12} & 1 & \multicolumn{2}{|l|}{0.0000109295947943995672548858} & 47,468 & 27,767 & 75,235 & 21.3164 & 2,615.57 & 52,148.93 & 16.48140056 \\ 
	\cline{2-11}		
	 & 2 & \multicolumn{2}{|l|}{0.0000036647139601286270917076} & 79,335 & 46,408 & 125,243 & 22.9227 & 4,367.86 & 87,158.41 & 18.05786797 \\ 
	\hline

	\multicolumn{11}{|c|}{} \\	
	\hline

	13 & 1 & & 0.00000360021965321270308169 & 111,202 & 65,049 & 176,251 & 23.2782 & 6,120.06 & 122,167.88 & 18.08348364 \\ 
	\hline

	\multicolumn{11}{|c|}{} \\	
	\hline

	14 & 1 & \multicolumn{2}{|l|}{0.000000064507500611466680552} & 190,537 & 111,457 & 301,994 & 27.8386 & 10,482.12 & 209,326.29 & 23.88595784 \\ 
	\hline

	\multicolumn{11}{|c|}{} \\	
	\hline

	\multirow{6}*{15} & 1 & & 0.0000035357119202803846457365 & 301,739 & 176,506 & 478,245 & 24.2944 & 16,596.18 & 331,494.17 & 18.10956784 \\ 
	\cline{2-11}
	 & 2 & & 0.0000034712041950929883664989 & 492,276 & 287,963 & 780,239 & 24.8023 & 27,072.05 & 540,820.46 & 18.13613234 \\  
	\cline{2-11}
	 & 3 & & 0.0000034066964668994453855464 & 682,813 & 399,233 & 1,082,233 & 25.1483 & 37,547.84 & 750,146.75 & 18.16319516  \\ 
	\cline{2-11}
	 & 4 & & 0.0000033421887464508240244483 & 873,350 & 510,877 & 1,384,227 & 25.4135 & 48,023.58 & 959,473.04 & 18.19077536 \\ 
	\cline{2-11}
	 & $\cdots$ & &  &  &  &  &  &  & & \\ 
	\cline{2-11}	 
	 & 55 & & 0.0000000523005213585975033236 & $\cdots$ & $\cdots$ & $\cdots$ & 32.0663 & 582,282 & 11,635,114 & 24.18859943 \\  
	\hline
					
\end{tabular}
}
\end{center}
\caption{Nodes - 3x+1 - Nodes 8 to 15}
\label{NodesProblem_3xPlus1_2}
\normalsize
\end{table}

\end{document}